\documentclass[pdflatex,sn-mathphys-num]{sn-jnl}

\usepackage{graphicx}%
\usepackage{multirow}%
\usepackage{amsmath,amssymb,amsfonts}%
\usepackage{amsthm}%
\usepackage{mathrsfs}%
\usepackage[title]{appendix}%
\usepackage{xcolor}%
\usepackage{textcomp}%
\usepackage{manyfoot}%
\usepackage{booktabs}%
\usepackage{algorithm}%
\usepackage{algorithmicx}%
\usepackage{algpseudocode}%
\usepackage{listings}%
\usepackage{array}
\usepackage{tikz}
\usetikzlibrary{
  arrows.meta,
  calc,
  positioning,
  decorations.pathreplacing,
  patterns,
  fit,
  backgrounds,
  intersections,
  shapes.geometric  
}

\hypersetup{
  colorlinks=true,
  linkcolor=blue!60!black,
  citecolor=blue!60!black,
  urlcolor=blue!60!black
}

\newtheorem{theorem}{Theorem}[section]
\newtheorem{proposition}[theorem]{Proposition}

\theoremstyle{definition}
\newtheorem{definition}[theorem]{Definition}
\newtheorem{remark}[theorem]{Remark}

\newtheorem{solution*}{Solution}[section]
\newcommand{\R}{\mathbb{R}}
\newcommand{\Gr}{\operatorname{Gr}}
\newcommand{\dom}{\operatorname{dom}}
\newcommand{\SOL}{\operatorname{SOL}}
\newcommand{\VI}{\operatorname{VI}}
\newcommand{\proj}{\Pi}

\newcommand{\argmin}{\operatorname*{argmin}}

\newcommand{\Fset}{\mathcal{F}}

\title{Introduction  to Mathematical Programming with Equilibrium Constraints (MPECs) and Bilevel Optimization}
\author{Louis Shuo Wang\thanks{Email: wang.s41@northeastern.edu}}
\affil{Department of Mathematics, Northeastern University, Boston, 02115, MA, USA}

\abstract{
Our aim is to explain mathematical programs with equilibrium constraints (MPECs), motivate them through applications, present the main equivalent formulations of equilibrium constraints, and summarize the basic existence theory for optimal solutions. The central message is that an MPEC is an optimization problem whose feasible set is partly defined by another optimization, variational inequality, complementarity system, or equilibrium model.}

\begin{document}
\maketitle

\section{Introduction and Big Picture}

An MPEC is an optimization problem in which part of the feasible set is defined implicitly through an equilibrium problem. The upper level chooses variables $(x,y)$ to optimize a global objective, while lower-level feasibility of $y$ is determined by a variational inequality, complementarity problem, or optimization problem parameterized by $x$.

This structure arises naturally in economics, engineering design, transportation, mechanics, robotics, and game theory. It generalizes bilevel programming and complementarity-constrained optimization. It also introduces major analytical challenges: nonconvex feasible regions, possible loss of closedness, nonsmoothness, and multi-valued lower-level response maps.

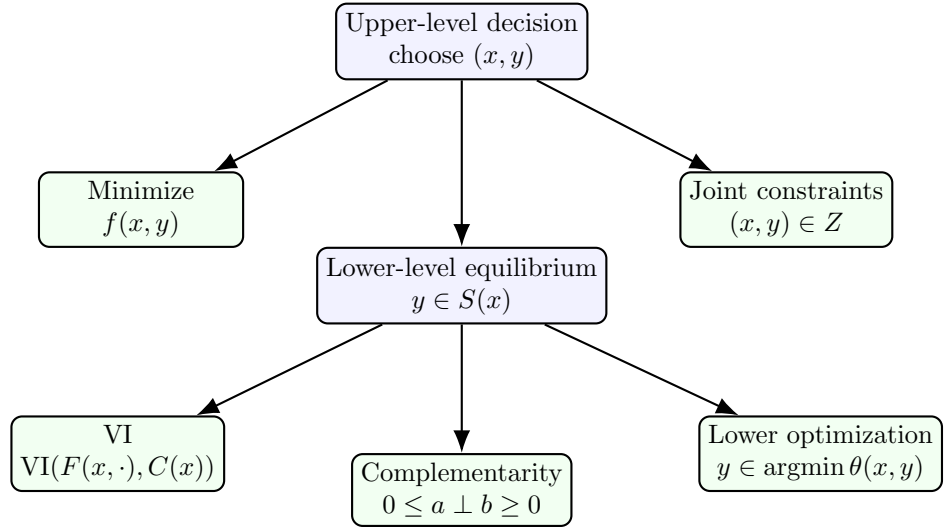
\begin{center}
\begin{tikzpicture}[
  node distance=1.7cm,
  block/.style={draw, rounded corners, thick, align=center, minimum width=3.3cm, minimum height=1.0cm, fill=blue!5},
  smallblock/.style={draw, rounded corners, thick, align=center, minimum width=2.7cm, minimum height=0.9cm, fill=green!5},
  arrow/.style={-{Latex[length=3mm]}, thick}
]
\node[block] (upper) {Upper-level decision\\ choose $(x,y)$};
\node[smallblock, below left=of upper] (objective) {Minimize\\ $f(x,y)$};
\node[smallblock, below right=of upper] (joint) {Joint constraints\\ $(x,y)\in Z$};
\node[block, below=2.2cm of upper] (lower) {Lower-level equilibrium\\ $y\in S(x)$};
\node[smallblock, below left=of lower] (vi) {VI\\ $\VI(F(x,\cdot),C(x))$};
\node[smallblock, below=of lower] (cp) {Complementarity\\ $0\le a\perp b\ge0$};
\node[smallblock, below right=of lower] (opt) {Lower optimization\\ $y\in\argmin\theta(x,y)$};
\draw[arrow] (upper) -- (objective);
\draw[arrow] (upper) -- (joint);
\draw[arrow] (upper) -- (lower);
\draw[arrow] (lower) -- (vi);
\draw[arrow] (lower) -- (cp);
\draw[arrow] (lower) -- (opt);
\end{tikzpicture}
\end{center}

\section{A One-Page Workshop Recap}

\begin{center}
\begin{tabular}{>{\raggedright\arraybackslash}p{0.25\textwidth}p{0.65\textwidth}}
\toprule
Topic & Key message \\
\midrule
MPEC definition & Optimize an upper-level objective subject to a lower-level equilibrium. \\
VI form & $y\in C(x)$ and $(v-y)^TF(x,y)\ge0$ for all $v\in C(x)$. \\
Complementarity form & Cone-constrained VIs become generalized complementarity systems. \\
KKT form & Valid when lower-level convex constraints admit multipliers under an appropriate CQ. \\
Fischer--Burmeister form & Use $\phi(a,b)=0$ to encode $0\le a\perp b\ge0$; for $g\le0$, use slack $u=-g\ge0$. \\
Implicit form & If $S(x)=\{y(x)\}$, the MPEC becomes an optimization problem in $x$ alone. \\
Existence & Need continuous objective, closed feasible region, and bounded lower level set. \\
Main warning & Closedness of $Z$ is not enough; the graph of the equilibrium map must also be closed. \\
\bottomrule
\end{tabular}
\end{center}

\section{Problem Formulation}

\subsection{The general MPEC model}
Let
\[
 f: \R^{n+m} \to \R, \qquad F: \R^{n+m} \to \R^m,
\]
let $Z \subset \R^{n+m}$ be nonempty and closed, and let
\[
 C: \R^n \rightrightarrows \R^m
\]
be a set-valued map with nonempty closed convex values on the domain of interest. Define
\[
 X = \{x \in \R^n : (x,y) \in Z \text{ for some } y \in \R^m\},
\]
and assume $X \subset \dom(C)$.

For each $x \in X$, let $S(x)$ denote the solution set of the variational inequality
\[
 y \in C(x), \qquad (v-y)^T F(x,y) \ge 0 \quad \forall v \in C(x).
\]
Then the MPEC is
\begin{equation}
\begin{aligned}
\min_{x,y} \quad & f(x,y)\\
\text{s.t.} \quad & (x,y) \in Z,\\
& y \in S(x).
\end{aligned}
\tag{MPEC}
\end{equation}
Its feasible region is
\[
 \Fset = Z \cap \Gr(S),
\]
where
\[
 \Gr(S)=\{(x,y)\in\R^{n+m}: y\in S(x)\}.
\]

\subsection{Interpretation of the model ingredients}
\begin{enumerate}[label=(\roman*)]
  \item $f(x,y)$ is the upper-level objective.
  \item $Z$ collects joint constraints on $(x,y)$.
  \item $F(x,y)$ defines the lower-level equilibrium condition.
  \item $C(x)$ is the lower-level feasible set for $y$ given $x$.
  \item $S(x)$ is the lower-level solution map, also called the reaction map.
\end{enumerate}


\subsection{Bilevel programming as a special case}
If there exists a continuously differentiable function $\theta(x,y)$ such that
\[
F(x,y)=\nabla_y \theta(x,y),
\]
then the lower-level VI is the first-order optimality condition for
\[
\min_y \; \theta(x,y) \quad \text{s.t. } y\in C(x).
\]
If $\theta(x,\cdot)$ is convex and an appropriate constraint qualification holds, the VI solution set coincides with the optimizer set of the lower-level problem.

\begin{remark}[Optimistic versus pessimistic bilevel modeling]
If the lower-level problem has multiple optimizers, the condition $y\in S(x)$ does not by itself specify how the lower-level solution is selected. The formulation above is the standard optimistic or strong version when the upper level is allowed to choose the most favorable $y\in S(x)$. A pessimistic or weak version instead evaluates the upper-level objective against the least favorable lower-level optimizer.
\end{remark}

\subsection{Complementarity systems as a special case}
If $C(x)$ is a convex cone, then the VI is equivalent to a generalized complementarity problem:
\[
 y\in C(x), \qquad F(x,y)\in C(x)^*, \qquad y^T F(x,y)=0,
\]
where the positive dual cone is
\[
 C(x)^* = \{z: z^T v \ge 0 \ \forall v\in C(x)\}.
\]
This sign convention matches the VI convention $(v-y)^TF(x,y)\ge0$. Under a polar-cone convention, the same statement is written with the corresponding sign change.

When $C(x)=\R^{m_1}\times \R_+^{m_2}$, one obtains mixed complementarity systems. When $m_1=0$, one obtains standard nonlinear complementarity problems.

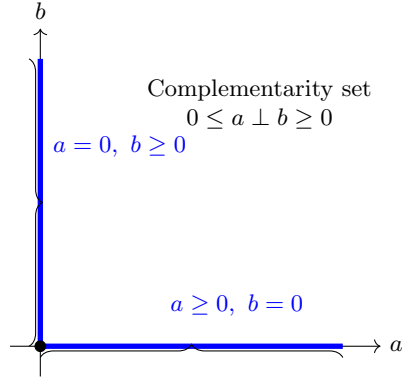
\begin{center}
\begin{tikzpicture}[scale=1.0, every node/.style={font=\small}]
\draw[->] (-0.4,0) -- (4.5,0) node[right] {$a$};
\draw[->] (0,-0.4) -- (0,4.2) node[above] {$b$};
\draw[line width=2pt, blue] (0,0) -- (4,0);
\draw[line width=2pt, blue] (0,0) -- (0,3.8);
\filldraw[black] (0,0) circle (2pt);
\node[blue, align=center] at (2.6,0.55) {$a\ge0,\ b=0$};
\node[blue, align=center] at (1.05,2.65) {$a=0,\ b\ge0$};
\node[align=center] at (2.9,3.2) {Complementarity set\\ $0\le a\perp b\ge0$};
\draw[decorate,decoration={brace,amplitude=5pt}] (0,-0.15) -- (4,-0.15);
\draw[decorate,decoration={brace,amplitude=5pt,mirror}] (-0.15,0) -- (-0.15,3.8);
\end{tikzpicture}
\end{center}

\subsection{Why MPECs are difficult}
The main difficulties are structural rather than cosmetic.
\begin{enumerate}[label=(\roman*)]
    \item \textbf{Nonconvexity.} Even if each lower-level problem is convex, the graph of the solution map $S$ can be highly nonconvex, so $\Fset$ is typically nonconvex.
    \item \textbf{Possible failure of closedness.} Although $Z$ may be closed, $\Gr(S)$ may fail to be closed unless the lower-level solution map has suitable closed-graph properties.
    \item \textbf{Nonsmoothness.} Even when $S(x)$ is single-valued, it may fail to be differentiable in $x$.
    \item \textbf{Multivaluedness.} The lower-level solution map may be set-valued, which creates ambiguity in the upper-level model and motivates optimistic and pessimistic Stackelberg interpretations.
    \item \textbf{Combinatorial structure.} Complementarity constraints encode disjunctions. For example, $0\le a\perp b\ge0$ means either $a=0$ or $b=0$ or both. This creates a latent finite-union structure similar to disjunctive programming.
\end{enumerate}

\subsection{Teaching summary}
\begin{center}
\begin{tabular}{>{\raggedright\arraybackslash}p{0.28\textwidth}p{0.62\textwidth}}
\toprule
Concept & Meaning \\
\midrule
Upper level & Global optimization over $(x,y)$ \\
Lower level & Equilibrium or optimization problem in $y$ parameterized by $x$ \\
Reaction map $S(x)$ & Set of equilibrium responses to $x$ \\
Feasible set $\Fset$ & Jointly feasible pairs that also satisfy lower-level equilibrium \\
Main difficulty & Nonconvexity, nonsmoothness, possible lack of closedness, and finite-union/combinatorial structure \\
\bottomrule
\end{tabular}
\end{center}

\section{Representative Problems}

This section motivates MPEC through representative applications. The common modeling template is
\[
\text{design, policy, or strategy choice} \quad + \quad \text{equilibrium response}.
\]

\subsection{Bracken--McGill bilevel programs}
These are early bilevel models from defense, production, and marketing. A firm or agency chooses upper-level decisions while accounting for adversarial or competitive lower-level choices. The lower-level variables represent the competitors' or opponents' optimizing reactions.

\subsection{Stackelberg games}
A leader chooses $x$, anticipating that followers choose an equilibrium response $y \in S(x)$. The leader solves
\[
\min_{x,y} f(x,y) \quad \text{s.t. } x\in X,\ y\in S(x).
\]
If $S(x)$ is multi-valued, two versions arise:
\begin{enumerate}[label=(\roman*)]
  \item \textbf{Strong/optimistic Stackelberg:} the leader assumes followers choose the response most favorable to the leader.
  \item \textbf{Weak/pessimistic Stackelberg:} the leader assumes followers choose the response least favorable to the leader.
\end{enumerate}

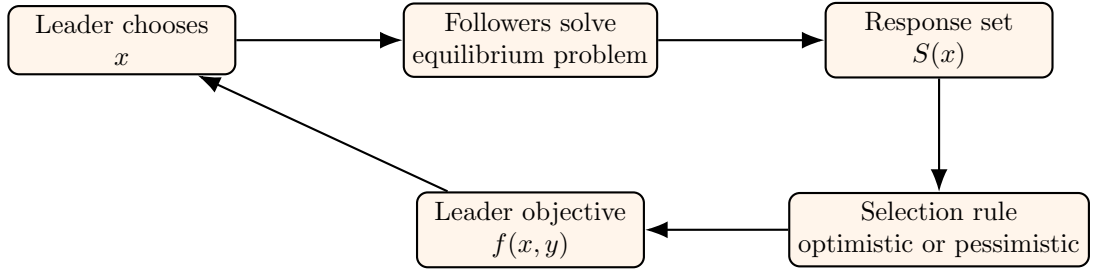
\begin{center}
\begin{tikzpicture}[
  node distance=1.7cm,
  box/.style={draw, rounded corners, thick, align=center, fill=orange!8, minimum width=3.0cm, minimum height=0.9cm},
  arrow/.style={-{Latex[length=3mm]}, thick}
]
\node[box] (leader) {Leader chooses\\ $x$};
\node[box, right=2.2cm of leader] (followers) {Followers solve\\ equilibrium problem};
\node[box, right=2.2cm of followers] (response) {Response set\\ $S(x)$};
\node[box, below=1.5cm of response] (selection) {Selection rule\\ optimistic or pessimistic};
\node[box, below=1.5cm of followers] (payoff) {Leader objective\\ $f(x,y)$};
\draw[arrow] (leader) -- (followers);
\draw[arrow] (followers) -- (response);
\draw[arrow] (response) -- (selection);
\draw[arrow] (selection) -- (payoff);
\draw[arrow] (payoff) -- (leader);
\end{tikzpicture}
\end{center}

\subsection{Misclassification minimization}
A minimum-misclassification hyperplane problem can be reformulated as a bilevel linear program by introducing lower-level variables that encode misclassification indicators through optimization.

\subsection{Robotics and contact mechanics}
In robotic motion with frictional contact, the dynamic constraints together with Coulomb friction and complementarity conditions lead to a complementarity system. If the original physical equations are infeasible or inconsistent, one may minimize a residual subject to the remaining equilibrium constraints. This produces an MPEC formulation.

\subsection{Residual minimization of complementarity systems}
Given an NCP,
\[
 x\ge 0, \qquad F(x)\ge 0, \qquad x^T F(x)=0,
\]
one may minimize a residual such as
\[
 \theta(x)=\frac12\|\min(x,F(x))\|^2.
\]
This residual minimization problem can itself be written as an MPEC or as a closely related nonsmooth optimization problem.

\subsection{Transportation and network design}
The upper level selects capacities, tolls, demands, or service frequencies. The lower level is a traffic equilibrium, route-choice equilibrium, or shortest-path problem. Representative examples include continuous network design, origin--destination demand adjustment, and transit planning.

\subsection{Facility location and production}
A large firm chooses plant locations and production quantities while anticipating market-clearing prices produced by a spatial price equilibrium problem. The lower-level equilibrium is often expressible as a complementarity system; the upper level optimizes profit.

\subsection{Structural mechanics and optimal design}
Design variables such as bar volumes or prestress forces appear at the upper level, while displacements, contact forces, and frictional forces are governed by complementarity systems at the lower level. These models are common in contact mechanics and shape optimization.

\section{Equivalent Constraint Formulations}

The VI constraint $y\in S(x)$ is compact but not always convenient. This section develops equivalent reformulations.

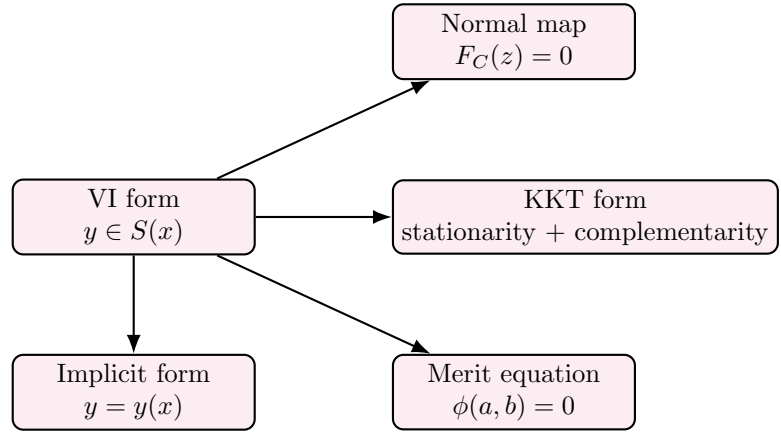
\begin{center}
\begin{tikzpicture}[
  node distance=1.3cm and 1.8cm,
  form/.style={draw, rounded corners, thick, align=center, fill=purple!7, minimum width=3.2cm, minimum height=0.9cm},
  arrow/.style={-{Latex[length=2.5mm]}, thick}
]
\node[form] (vi) {VI form\\ $y\in S(x)$};
\node[form, above right=of vi] (normal) {Normal map\\ $F_C(z)=0$};
\node[form, right=of vi] (kkt) {KKT form\\ stationarity + complementarity};
\node[form, below right=of vi] (merit) {Merit equation\\ $\phi(a,b)=0$};
\node[form, below=of vi] (implicit) {Implicit form\\ $y=y(x)$};
\draw[arrow] (vi) -- (normal);
\draw[arrow] (vi) -- (kkt);
\draw[arrow] (vi) -- (merit);
\draw[arrow] (vi) -- (implicit);
\end{tikzpicture}
\end{center}

\subsection{Normal-map formulation}
Let $C\subset \R^m$ be nonempty, closed, and convex. Denote by $\proj_C(z)$ the Euclidean projection of $z$ onto $C$. For a VI$(F,C)$, the associated normal map is
\[
F_C(z)=F(\proj_C(z))+z-\proj_C(z).
\]
Then
\[
 y \in \SOL(F,C) \iff \exists z \text{ such that } y=\proj_C(z),\ F_C(z)=0.
\]
Hence, for fixed $x$, the lower-level VI can be rewritten as an equation in a new variable $v$:
\[
F_{C(x)}(v)=F(x,\proj_{C(x)}(v))+v-\proj_{C(x)}(v)=0.
\]
This converts the lower-level VI into an equation, but introduces nonsmoothness through the projection operator.

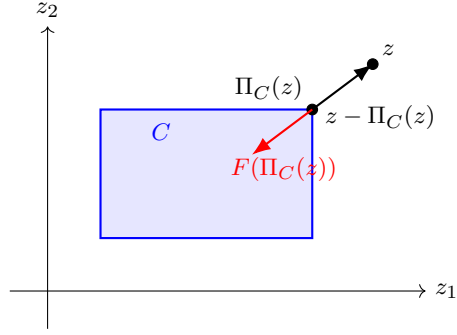
\begin{tikzpicture}[scale=1.0, every node/.style={font=\small}]
\draw[->] (-0.5,0) -- (5,0) node[right] {$z_1$};
\draw[->] (0,-0.5) -- (0,3.5) node[above] {$z_2$};

\fill[blue!10] (0.7,0.7) rectangle (3.5,2.4);
\draw[thick, blue] (0.7,0.7) rectangle (3.5,2.4);
\node[blue] at (1.5,2.1) {$C$};

\coordinate (z) at (4.3,3.0);
\coordinate (p) at (3.5,2.4);
\coordinate (f) at (2.7,1.8); 

\filldraw[black] (z) circle (2pt) node[above right] {$z$};
\filldraw[black] (p) circle (2pt) node[above left] {$\Pi_C(z)$};

\draw[thick, -{Latex[length=2.5mm]}] (p) -- (z);
\node at (4.4,2.3) {$z-\Pi_C(z)$};

\draw[thick, red, -{Latex[length=2.5mm]}] (p) -- (f);
\node[red, below right] at (2.3,1.9) {$F(\Pi_C(z))$};
\end{tikzpicture}

\subsection{KKT formulation}
Assume now that
\[
C(x)=\{y\in\R^m : g(x,y)\le 0\},
\]
where each $g_i(x,\cdot)$ is convex and differentiable. Under a suitable multiplier constraint qualification, such as Slater's condition for the lower-level convex feasible set or another condition guaranteeing KKT multipliers, the VI is equivalent to a KKT system: there exists $\lambda\in\R^\ell$ such that
\begin{align*}
F(x,y)+\sum_{i=1}^\ell \lambda_i \nabla_y g_i(x,y)&=0,\\
\lambda_i &\ge 0,\\
g_i(x,y)&\le 0,\\
\lambda_i g_i(x,y)&=0 \qquad (i=1,\dots,\ell).
\end{align*}
Thus the MPEC becomes
\begin{equation}
\begin{aligned}
\min_{x,y,\lambda} \quad & f(x,y)\\
\text{s.t.} \quad & (x,y)\in Z,\\
& F(x,y)+\sum_{i=1}^\ell \lambda_i \nabla_y g_i(x,y)=0,\\
& \lambda\ge 0,\quad g(x,y)\le 0,\quad \lambda_i g_i(x,y)=0\quad i=1,\dots,\ell.
\end{aligned}
\tag{KKT-MPEC}
\end{equation}
This is one of the most important equivalent forms in theory and computation.

\subsection{SBCQ and closedness of the feasible region}
A central issue is whether KKT multipliers can be chosen in a uniformly bounded way along feasible sequences. A sequentially bounded constraint qualification (SBCQ) says, roughly, that along any convergent feasible sequence, one can choose bounded multipliers.

Under continuity, convexity, and suitable multiplier-existence assumptions, SBCQ implies:
\begin{enumerate}[label=(\roman*)]
  \item $y\in S(x)$ if and only if there exists a KKT multiplier at $(x,y)$,
  \item the MPEC feasible region is closed.
\end{enumerate}
This is crucial for existence theory.

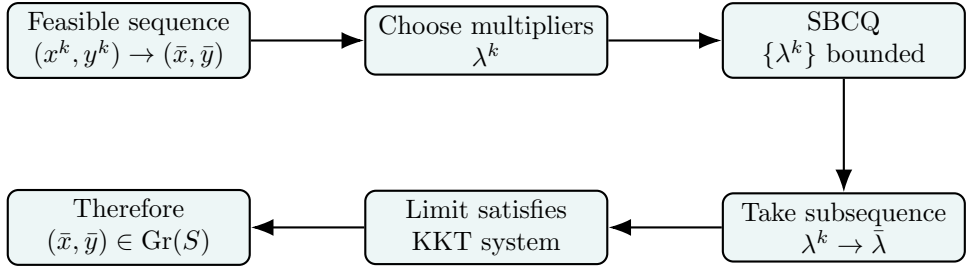
\begin{center}
\begin{tikzpicture}[
  node distance=1.5cm,
  box/.style={draw, rounded corners, thick, align=center, fill=teal!7, minimum width=3.2cm, minimum height=0.9cm},
  arrow/.style={-{Latex[length=3mm]}, thick}
]
\node[box] (seq) {Feasible sequence\\ $(x^k,y^k)\to(\bar x,\bar y)$};
\node[box, right=of seq] (mult) {Choose multipliers\\ $\lambda^k$};
\node[box, right=of mult] (bounded) {SBCQ\\ $\{\lambda^k\}$ bounded};
\node[box, below=of bounded] (subseq) {Take subsequence\\ $\lambda^k\to\bar\lambda$};
\node[box, left=of subseq] (limit) {Limit satisfies\\ KKT system};
\node[box, left=of limit] (closed) {Therefore\\ $(\bar x,\bar y)\in\Gr(S)$};
\draw[arrow] (seq) -- (mult);
\draw[arrow] (mult) -- (bounded);
\draw[arrow] (bounded) -- (subseq);
\draw[arrow] (subseq) -- (limit);
\draw[arrow] (limit) -- (closed);
\end{tikzpicture}
\end{center}

\subsection{Piecewise structure of the feasible set}
The complementarity conditions partition the feasible region into finitely many pieces according to which constraints are active and which multipliers are positive. Therefore, under suitable assumptions, the feasible set is a finite union of closed sets. This reveals the latent combinatorial structure of MPEC.

For each complementarity pair $0\le \lambda_i \perp u_i\ge0$, where $u_i=-g_i(x,y)$, the feasible alternatives are
\[
\lambda_i=0,\ u_i\ge0
\qquad\text{or}\qquad
\lambda_i\ge0,\ u_i=0.
\]
With $\ell$ complementarity pairs, this creates up to $2^\ell$ active-set regimes.

\subsection{Affine variational inequalities and MPAECs}
If
\[
F(x,y)=q+Nx+My, \qquad g(x,y)=b+Ax+By,
\]
then the lower-level VI is affine and the resulting problem is often called an MPEC with affine equilibrium constraints, or MPAEC. In this case, each active-set piece of the feasible set is polyhedral, so the global feasible set is a finite union of polyhedra.

This special case is much more structured and often serves as a first serious testbed for theory and algorithms.

\subsection{Merit-function reformulations}
Complementarity conditions can be converted into equations using an NCP function. A standard example is the Fischer--Burmeister function:
\[
\phi(a,b)=\sqrt{a^2+b^2}-(a+b).
\]
It satisfies
\[
\phi(a,b)=0 \iff a\ge 0,\ b\ge 0,\ ab=0.
\]
For the KKT inequalities $g_i(x,y)\le0$, define the nonnegative slack
\[
u_i(x,y):=-g_i(x,y)\ge0.
\]
Then complementarity is equivalently written as
\[
\phi(\lambda_i,u_i(x,y))=0,\qquad i=1,\dots,\ell.
\]
This yields an equation-based reformulation of MPEC and motivates smooth or semismooth numerical methods.

\subsection{Implicit-program formulation}
If, for each $x$, the lower-level VI has a unique solution $y(x)$, then the MPEC reduces to
\[
\min_x f(x,y(x)) \quad \text{s.t. } x\in X,\ (x,y(x))\in Z.
\]
A common sufficient condition for uniqueness is strict or strong monotonicity of $F(x,\cdot)$ on $C(x)$.

\begin{definition}[Monotonicity]
A mapping $G:K\to\R^m$ is:
\begin{enumerate}[label=(\roman*)]
  \item \emph{monotone} if $(u-v)^T(G(u)-G(v))\ge 0$ for all $u,v\in K$,
  \item \emph{strictly monotone} if the inequality is strict for $u\neq v$,
  \item \emph{strongly monotone} if there exists $c>0$ such that
  \[
  (u-v)^T(G(u)-G(v))\ge c\|u-v\|^2 \quad \forall u,v\in K.
  \]
\end{enumerate}
\end{definition}

\begin{proposition}[Existence and uniqueness for a VI]
Let $K\subset\R^m$ be nonempty, closed, and convex, and let $G:K\to\R^m$ be continuous. If $K$ is compact, then $\VI(G,K)$ has a solution. If, in addition, $G$ is strictly monotone, the solution is unique. If $K$ is unbounded, existence requires an additional coercivity or growth condition; strong monotonicity is a common condition that yields both uniqueness and, under the standard VI hypotheses, existence.
\end{proposition}

For complementarity problems on orthants, a related uniqueness condition is the uniform $P$ condition.

\subsection{Constraint qualifications that imply SBCQ}
SBCQ is implied by more familiar nonlinear programming conditions, including:
\begin{enumerate}[label=(\roman*)]
  \item Mangasarian--Fromovitz constraint qualification (MFCQ),
  \item constant-rank constraint qualification (CRCQ),
  \item linear independence constraint qualification (LICQ).
\end{enumerate}
LICQ is strongest among these and yields uniqueness of the multiplier vector.

\subsection{Teaching summary}
The most important takeaway is:
\begin{quote}
The same equilibrium constraint can be expressed in projection form, KKT form, complementarity form, merit-function form, or implicit form. Different forms are useful for different proofs and algorithms.
\end{quote}

\section{Existence of Optimal Solutions}

\subsection{General principle}
Existence of a minimizer usually follows from two ingredients:
\begin{enumerate}[leftmargin=1.5em]
  \item the feasible set is closed,
  \item at least one lower level set of the objective is nonempty and bounded.
\end{enumerate}
If these hold and the objective is continuous, a minimizer exists by the Weierstrass theorem.

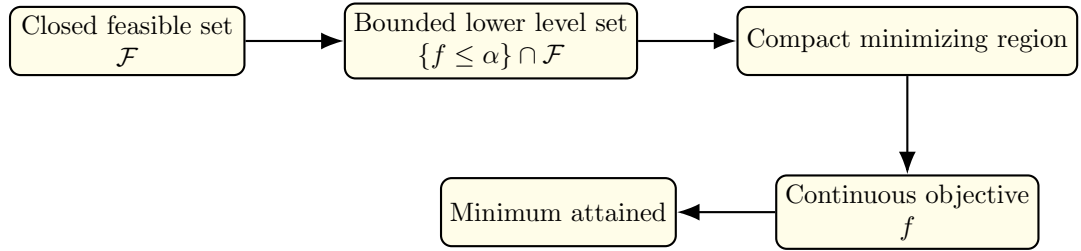
\begin{center}
\begin{tikzpicture}[
  node distance=1.3cm,
  box/.style={draw, rounded corners, thick, align=center, fill=yellow!10, minimum width=3.1cm, minimum height=0.9cm},
  arrow/.style={-{Latex[length=3mm]}, thick}
]
\node[box] (closed) {Closed feasible set\\ $\Fset$};
\node[box, right=of closed] (bounded) {Bounded lower level set\\ $\{f\le \alpha\}\cap\Fset$};
\node[box, right=of bounded] (compact) {Compact minimizing region};
\node[box, below=of compact] (cont) {Continuous objective\\ $f$};
\node[box, left=of cont] (attain) {Minimum attained};
\draw[arrow] (closed) -- (bounded);
\draw[arrow] (bounded) -- (compact);
\draw[arrow] (compact) -- (cont);
\draw[arrow] (cont) -- (attain);
\end{tikzpicture}
\end{center}

\subsection{Coercivity and bounded level sets}
A function $\psi$ on a set $S$ is coercive if
\[
\|z\|\to\infty,\ z\in S \implies \psi(z)\to\infty.
\]
A coercive continuous objective on a closed feasible set attains its minimum. More generally, it is enough that one lower level set is nonempty and bounded.

\subsection{Existence via the KKT formulation}
For the KKT reformulation, if the feasible set in $(x,y,\lambda)$-space is closed and one lower level set is nonempty and bounded, then an optimal solution exists.

For the original MPEC, a clean sufficient condition is the following.

\begin{theorem}[Basic existence result for MPEC]
Assume that:
\begin{enumerate}[label=(\roman*)]
  \item $f$ is continuous,
  \item $Z$ is closed,
  \item the lower-level data are continuous,
  \item each $g_i(x,\cdot)$ is convex and differentiable,
  \item a multiplier constraint qualification and SBCQ hold so that $\Gr(S)$ is closed,
  \item some lower level set $\{(x,y)\in\Fset:f(x,y)\le \alpha\}$ is nonempty and bounded.
\end{enumerate}
Then the MPEC has a global optimal solution.
\end{theorem}

\subsection{The affine/quadratic case}
In the MPAEC setting with polyhedral $Z$ and quadratic objective $f$, each active-set piece is polyhedral. An analog of the Frank--Wolfe theorem says that a quadratic objective over a finite union of polyhedra attains a global solution when it is bounded below on the feasible set, equivalently when it is bounded below on every feasible ray of every polyhedral piece.

\subsection{Interpretation}
The hard part of existence theory for MPEC is usually not continuity of the objective. The hard part is proving that the feasible region is closed despite being generated by equilibrium constraints. Once closedness and boundedness are available, existence follows from standard compactness arguments.

\section{Exercises}

\subsection{True or False}
\begin{enumerate}[label=(\roman*)]
    \item \textbf{Topology of the Feasible Region:}
    
    For a general Mathematical Program with Equilibrium Constraints (MPEC), the feasible region is guaranteed to be a closed set, even though it is typically nonconvex and potentially disconnected.
    \item \textbf{Complexity of Bilevel Programs:}

    A bilevel linear program—where the overall objective is linear, the joint constraints are polyhedral, and the inner problem is a linear program—belongs to the class of strongly NP-hard problems.
    \item \textbf{Stackelberg Game Formulations}

    In a Stackelberg game where the followers' rational reaction set $S(x)$ is multi-valued, the "strong" (optimistic) and "weak" (pessimistic) formulations yield equivalent objective values for the leader.
    \item \textbf{The Normal Map Reformulation:}

    Reformulating the inner variational inequality of an MPEC using the normal map equation (involving the projection operator $\Pi$) successfully transforms the problem into a standard, continuously differentiable equality-constrained optimization problem.
    \item \textbf{Constraint Qualifications:}

    The Sequentially Bounded Constraint Qualification (SBCQ) is a critical assumption because it guarantees that the feasible region of the MPEC is closed, provided the lower-level constraints are continuous and convex.
    \item \textbf{Constraint Qualification Hierarchy:}

    If the Mangasarian-Fromovitz Constraint Qualification (MFCQ) holds at all pairs in the joint feasible region, it implies that the SBCQ also holds.
    \item \textbf{Monotonicity and Implicit Forms:}

    If the equilibrium function $F(x, \cdot)$ is strongly monotone on a closed convex set $C(x)$ for each $x$, the reaction map $S(x)$ becomes single-valued, allowing the MPEC to be reformulated as an implicit optimization problem solely in terms of the upper-level variable $x$.
    \item \textbf{Existence of Optimal Solutions:}

    The residual minimization problem of any Affine Variational Inequality (AVI) over a nonempty polyhedron is guaranteed to have an optimal solution.
\end{enumerate}

\subsection{Computational questions}
\subsubsection*{Question 1: Topology of the Rational Set and SBCQ}
Consider the following bilevel program in $\mathbb{R}^2$:
\begin{align*}
    \text{minimize} \quad & f(x,y) \\
    \text{subject to} \quad & x \in [-1, 1], \\
    \text{and} \quad & y \in \text{argmin} \{ v \mid -1 \le v \le 1, \; xv \le 0 \}.
\end{align*}

\begin{enumerate}
    \item[(a)] Compute the exact analytic expression for the reaction map $S(x)$ for all $x \in [-1, 1]$.
    \item[(b)] Graph or explicitly define the joint feasible region $\mathcal{F} = \mathcal{Z} \cap \text{Gr}(S)$. Prove mathematically whether $\mathcal{F}$ is closed and/or convex.
    \item[(c)] The Sequentially Bounded Constraint Qualification (SBCQ) is critical for guaranteeing the closedness of $\mathcal{F}$. Evaluate the KKT multipliers for the inner optimization problem at the sequence of points $x_k = -\frac{1}{k}$ as $k \to \infty$. Prove that the SBCQ fails at $x=0$, and explain how this computational failure relates to your topological findings in part (b).
\end{enumerate}

\subsubsection*{Question 2: The Normal Map Reformulation}
You are tasked with computationally solving a simple MPEC where the inner problem is a nonlinear complementarity problem (NCP) over the non-negative orthant $\mathbb{R}^m_+$. 
\begin{align*}
    \text{minimize} \quad & \frac{1}{2} \| x - y \|^2 \\
    \text{subject to} \quad & x \ge 0, \\
    \text{and} \quad & y \ge 0, \quad F(x,y) \ge 0, \quad y^T F(x,y) = 0.
\end{align*}
Assume $F(x,y) = M y + N x + q$, where $M \in \mathbb{R}^{m \times m}$ is a positive definite matrix.

\begin{enumerate}
    \item[(a)] Define the projection operator $\Pi_{\mathbb{R}^m_+}(z)$ and construct the normal map $F_{\mathbb{R}^m_+}(z)$ associated with the lower-level variational inequality. 
    \item[(b)] Perform a change of variables to reformulate the above MPEC into its equivalent \textbf{Normal Form} (an optimization problem with a single non-smooth equality constraint). Explicitly state the new objective function and constraints in terms of your new variable $z$.
    \item[(c)] Discuss the primary computational hurdle of using standard Newton-type methods on your formulation in (b), and briefly suggest how one might overcome the non-differentiability of the projection operator in a numerical solver.
\end{enumerate}

\subsubsection*{Question 3: Constraint Qualifications (MFCQ vs. CRCQ)}
Let $\mathcal{Z} = \mathbb{R}^4$ and let the lower-level constraint set $C(x)$ be defined by the inequalities $g_1(x,y) \le 0$ and $g_2(x,y) \le 0$, where for $x \in \mathbb{R}^2, y \in \mathbb{R}^2$:
\begin{align*}
    g_1(x,y) &= y_1 + (y_2)^2 - x_1 \\
    g_2(x,y) &= y_1 - x_2
\end{align*}

\begin{enumerate}
    \item[(a)] Formulate the Jacobian matrices $\nabla_y g_1(x,y)$ and $\nabla_y g_2(x,y)$.
    \item[(b)] Evaluate the linear independence of the gradients at the point $(x,y) = (0,0,0,0)$. 
    \item[(c)] Prove that the Mangasarian-Fromovitz Constraint Qualification (MFCQ) holds at $(0,0,0,0)$, but the Constant-Rank Constraint Qualification (CRCQ) fails. Show your matrix rank computations.
\end{enumerate}

\subsubsection*{Question 4: Monotonicity and the Implicit Form (25 points)}
Consider a transportation continuous network design problem formulated as an MPEC, where the lower-level problem determines link flows $y \in \mathbb{R}^m$ given capacity enhancements $x \in \mathbb{R}^n$. The lower-level problem is an Affine Variational Inequality (AVI) defined by $F(x,y) = M y + N x + q$ over a polyhedral set $C(x) = \{y \ge 0 \mid A y \le b(x)\}$.

\begin{enumerate}
    \item[(a)] Define what it means for $F(x, \cdot)$ to be \textit{strongly monotone} on $C(x)$. State the necessary and sufficient algebraic condition on the matrix $M$ for this to hold.
    \item[(b)] If $M$ satisfies the condition in (a), what topological property does the solution set $S(x)$ possess? 
    \item[(c)] Assuming strong monotonicity holds and $C(x)$ is independent of $x$ (i.e., $C(x) = Y$), write down the equivalent \textbf{Implicit Form} of the MPEC. Why is this form computationally advantageous, and what is its main limitation if joint constraints $(x,y) \in \mathcal{Z}$ are introduced?
\end{enumerate}

\subsection{Solutions}
\subsubsection*{Solutions for True or False}
\begin{enumerate}[label=(\roman*)]
    \item False. Example 1.1.2 explicitly demonstrates that the feasible region of an MPEC can be disconnected and not closed. The lack of closedness is a primary reason why MPECs are so intractable without assuming constraint qualifications like the SBCQ.
    \item True. As stated in Section 1.1, despite its apparent simplicity (linear objective, polyhedral constraints, linear inner problem), the bilevel linear program is strongly NP-hard, meaning no fully polynomial approximation scheme exists unless P = NP.
    \item False. When S(x) is a singleton, the two are identical. However, when S(x) is multi-valued, the strong formulation assumes followers choose the reaction most favorable to the leader (minimizing the leader's cost over S(x)), while the weak formulation assumes they choose the least favorable reaction (taking the supremum over S(x)). These are generally not equivalent.
    \item False. While the normal map formulation does transform the VI into a single equation, the projection operator used in this formulation is generally not Fréchet or Gâteaux differentiable. Therefore, it results in a non-smooth/nondifferentiable optimization problem.
    \item True. Theorem 1.3.4 explicitly proves that under the SBCQ (and assuming convexity/continuity of the lower-level constraint functions), the feasible region of the MPEC is a closed set. SBCQ ensures that bounded KKT multipliers exist for convergent sequences.
    \item True. Proposition 1.3.6 establishes that MFCQ is a stronger condition that implies the SBCQ. Roughly speaking, MFCQ ensures that all multipliers are bounded on bounded sets, which easily satisfies the sequential boundedness required by SBCQ.
    \item True. Strong monotonicity guarantees that the inner variational inequality has exactly one solution (Proposition 1.3.10c). Because S(x) is a singleton, let's call it y(x), we can substitute it into the upper-level objective to get f(x,y(x)), reducing the problem to an implicit program in x.
    \item True. Corollary 1.4.4 proves this. By reformulating the AVI residual minimization into an MPAEC (Mathematical Program with Affine Equilibrium Constraints) with a convex quadratic objective, the problem is bounded below by zero. Extending the Frank-Wolfe theorem (Theorem 1.4.3), a quadratic function bounded below on the feasible rays of a polyhedral set attains its minimum.
\end{enumerate}

\subsubsection*{Solutions for computational questions}
\paragraph{Question 1: Topology of the Rational Set and SBCQ (25 points)}

\textbf{Part (a): Compute $S(x)$ (8 points)} \\
The inner problem is: $\min y$ subject to $y \in [-1, 1]$ and $xy \le 0$. 
\begin{itemize}
    \item If $x \in (0, 1]$, $xy \le 0 \implies y \le 0$. The feasible set for $y$ is $[-1, 0]$. The minimum is at $y = -1$.
    \item If $x \in [-1, 0)$, $xy \le 0 \implies y \ge 0$. The feasible set for $y$ is $[0, 1]$. The minimum is at $y = 0$.
    \item If $x = 0$, $xy \le 0 \implies 0 \le 0$ (always true). The feasible set for $y$ is $[-1, 1]$. The minimum is at $y = -1$.
\end{itemize}
Therefore, the reaction map is:
$$ S(x) = \begin{cases} \{-1\} & \text{if } x \in [0, 1] \\ \{0\} & \text{if } x \in [-1, 0) \end{cases} $$
\textit{Rubric: 4 pts for identifying the correct intervals. 4 pts for the correct $y$ values.}

\textbf{Part (b): Feasible Region $\mathcal{F}$ (8 points)} \\
The joint feasible region is $\mathcal{F} = \{ (x,-1) \mid x \in [0, 1] \} \cup \{ (x,0) \mid x \in [-1, 0) \}$. \\
\textbf{Closedness:} $\mathcal{F}$ is \textbf{not closed}. Take the sequence $x_k = -1/k$. As $k \to \infty$, $x_k \to 0$ and $y_k \to 0$. The limit point is $(0,0)$. However, $S(0) = \{-1\}$, so $(0,0) \notin \mathcal{F}$. \\
\textbf{Convexity:} $\mathcal{F}$ is \textbf{not convex}. It is the union of two disjoint, non-collinear line segments. \\
\textit{Rubric: 2 pts for the correct set definition. 3 pts for proving it is not closed (limit point argument). 3 pts for proving it is not convex.}

\textbf{Part (c): SBCQ Failure (9 points)} \\
The Lagrangian for the inner problem is $L(y, \lambda, \mu) = y + \lambda(xy) + \mu_1(y-1) + \mu_2(-y-1)$. \\
The stationarity KKT condition is $\nabla_y L = 1 + \lambda x + \mu_1 - \mu_2 = 0$. \\
At $x_k = -1/k$ and $y_k = 0$, the constraints $y \le 1$ and $y \ge -1$ are strictly inactive, so $\mu_1 = \mu_2 = 0$. The constraint $xy \le 0$ is active ($0 \le 0$). 
Substitute into the KKT condition:
$$ 1 + \lambda_k \left(-\frac{1}{k}\right) = 0 \implies \lambda_k = k $$
As $k \to \infty$, $\lambda_k \to \infty$. Thus, the multipliers are unbounded for a sequence converging to a feasible limit point, meaning \textbf{SBCQ fails} at $x=0$. This computational failure perfectly aligns with the topological failure in (b), as SBCQ is the sufficient condition for $\mathcal{F}$ to be closed. \\
\textit{Rubric: 3 pts for correct KKT setup. 4 pts for calculating $\lambda_k = k$. 2 pts for linking the unbounded multiplier to the lack of closedness.}

\paragraph{Question 2: The Normal Map Reformulation (25 points)}

\textbf{Part (a): Operator and Normal Map (8 points)} \\
The projection operator onto the non-negative orthant is the component-wise maximum: $\Pi_{\mathbb{R}^m_+}(z) = \max(z, 0) = z^+$. \\
The normal map associated with the pair $(F, \mathbb{R}^m_+)$ is defined as $F_{\mathbb{R}^m_+}(z) = F(z^+) - z^-$, where $z^- = \max(-z, 0)$. \\
Substituting $F(x,y)$:
$$ F_{\mathbb{R}^m_+}(x, z) = M z^+ + N x + q - z^- $$
\textit{Rubric: 4 pts for the correct projection operator. 4 pts for the correct normal map expression.}

\textbf{Part (b): Reformulated MPEC (9 points)} \\
Using the substitution $y = z^+$, the new optimization problem is over variables $(x,z)$:
\begin{align*}
    \text{minimize} \quad & \frac{1}{2} \| x - z^+ \|^2 \\
    \text{subject to} \quad & x \ge 0, \\
    \text{and} \quad & M z^+ + N x + q - z^- = 0.
\end{align*}
\textit{Rubric: 4 pts for substituting $y=z^+$ in the objective. 5 pts for accurately replacing the complementarity constraints with the single normal map equation.}

\textbf{Part (c): Computational Hurdles (8 points)} \\
The primary hurdle is that the operators $z^+$ and $z^-$ are \textbf{non-differentiable} at $z_i = 0$. Consequently, the normal map equation is non-smooth, meaning standard gradient-based or Newton methods cannot be directly applied. \\
\textbf{Overcoming this:} One can use \textit{Semismooth Newton methods} (which rely on generalized Jacobians/subdifferentials), or apply a \textit{smoothing function} (e.g., replacing $\max(0,z_i)$ with a smooth approximation like the Fischer-Burmeister function or a softplus function) to restore differentiability. \\
\textit{Rubric: 4 pts for identifying non-differentiability at $z=0$. 4 pts for a valid algorithmic solution (semismooth Newton, smoothing functions, etc.).}

\paragraph{Question 3: Constraint Qualifications (25 points)}

\textbf{Part (a): Jacobians (6 points)} \\
$$ \nabla_y g_1(x,y) = \begin{pmatrix} 1 \\ 2y_2 \end{pmatrix}, \quad \nabla_y g_2(x,y) = \begin{pmatrix} 1 \\ 0 \end{pmatrix} $$
\textit{Rubric: 3 pts for each correct gradient vector.}

\textbf{Part (b): Linear Independence at the Origin (6 points)} \\
At $(x,y) = (0,0,0,0)$, we have:
$$ \nabla_y g_1(0,0) = \begin{pmatrix} 1 \\ 0 \end{pmatrix}, \quad \nabla_y g_2(0,0) = \begin{pmatrix} 1 \\ 0 \end{pmatrix} $$
These two vectors are identical and therefore \textbf{linearly dependent}. \\
\textit{Rubric: 3 pts for evaluating the gradients at the origin. 3 pts for stating they are linearly dependent.}

\textbf{Part (c): MFCQ vs CRCQ (13 points)} \\
\textbf{MFCQ:} MFCQ holds if there exists a vector $v \in \mathbb{R}^2$ such that $v^T \nabla_y g_1 < 0$ and $v^T \nabla_y g_2 < 0$. 
Let $v = \begin{pmatrix} -1 \\ 0 \end{pmatrix}$. 
$$ v^T \begin{pmatrix} 1 \\ 0 \end{pmatrix} = -1 < 0 $$
Thus, MFCQ holds at the origin. \\
\textbf{CRCQ:} CRCQ requires the rank of the set of active constraint gradients to remain constant in a neighborhood around the point.
The Jacobian matrix of the active constraints is $J = \begin{pmatrix} 1 & 1 \\ 2y_2 & 0 \end{pmatrix}$. 
At the origin ($y_2 = 0$), the rank of $J$ is 1. However, for any arbitrarily small neighborhood where $y_2 \neq 0$, the determinant is $-2y_2 \neq 0$, meaning the rank jumps to 2. Because the rank is not constant in the neighborhood, \textbf{CRCQ fails}. \\
\textit{Rubric: 6 pts for proving MFCQ (finding a valid vector $v$). 7 pts for proving CRCQ fails (showing the rank changes from 1 to 2 in the neighborhood).}

\paragraph{Problem 4: Monotonicity and the Implicit Form}

\textbf{Part (a): Strong Monotonicity (8 points)} \\
$F(x, \cdot)$ is strongly monotone on $C(x)$ if there exists a constant $c > 0$ such that for all $u, v \in C(x)$:
$$ (u - v)^T (F(x,u) - F(x,v)) \ge c \|u - v\|^2 $$
Since $F$ is affine with respect to $y$, $F(x,u) - F(x,v) = M(u-v)$. The condition becomes $(u-v)^T M (u-v) \ge c \|u-v\|^2$. This holds if and only if $M$ is a \textbf{positive definite matrix}. \\
\textit{Rubric: 4 pts for the mathematical definition of strong monotonicity. 4 pts for stating $M$ must be positive definite.}

\textbf{Part (b): Topological Property of $S(x)$ (8 points)} \\
If strong monotonicity holds, the solution set $S(x)$ consists of \textbf{exactly one element} (it is a singleton). The reaction map is single-valued. \\
\textit{Rubric: 8 pts for stating the solution is unique/a singleton.}

\textbf{Part (c): The Implicit Form (9 points)} \\
Because $S(x)$ is single-valued, let $y(x)$ denote this unique solution vector. The implicit form of the MPEC is:
\begin{align*}
    \text{minimize} \quad & f(x, y(x)) \\
    \text{subject to} \quad & x \in X
\end{align*}
\textbf{Advantage:} It reduces the MPEC to a standard nonlinear program (NLP) in the variable $x$ alone, stripping away the combinatorial complexity of the complementarity constraints. \\
\textbf{Limitation:} If there are joint constraints $(x,y) \in \mathcal{Z}$ (rather than just $x \in X$), the implicit function $y(x)$ might violate these constraints. It is extremely difficult to characterize the feasible region of $x$ such that $(x, y(x)) \in \mathcal{Z}$, which can lead to an empty feasible set or highly intractable constraints. \\
\textit{Rubric: 3 pts for the correct formulation. 3 pts for the advantage (reduces to NLP). 3 pts for the limitation (handling joint constraints).}

\section{Bibliographic remarks and Acknowledgment}

\textbf{This note is mainly based on Chapter 1, in the MPEC monograph of Zhi-Quan Luo, Jong-Shi Pang and Daniel Ralph.} Please see the relevant references: 
\cite{falk1995bilevel,
aghasi2025fully,
hong2023two,
kovccvara1994optimization,
chaudet2020shape,
liu2023auxiliary,
liu2025bidirectional,
outrata1995numerical,
cui2023complexity,
christof2020nonsmooth,
rawat2026augmented,
robinson1980strongly,
chen2025control,
shin2023near,
bolte2024differentiating,
nghia2025geometric,
wang2025analysis,
chen2025aubin,
liu2024auxiliary,
kojima1980strongly,
chen2026characterizations,
cui2026lipschitz,
shin2022exponential,
de2023function,
bank1982d,
wang2026algebraic,
bonnans1994local,
khanh2024globally,
chen2025two,
mohammadi2022variational,
dussault2026polyhedral,
gowda1994stability,
liu2025learning,
qi2000constant,
facchinei1998accurate,
jittorntrum2009solution,
kyparisis1992parametric,
liu1995sensitivity,
pang11995stability,
liu2025risk,
liu2024feasibility,
qiu1992sensitivity,
aussel2024variational,
reinoza1985strong,
facchinei2003finite,
harker1990finite,
kyparisis1990sensitivity,
wang2025multi,
kleinmichel1972av,
ortega2000iterative,
wang2025analysis1,
bai2021matrix,
lin2026hierarchical,
gao2022rolling,
gander2026landmarks,
mishchenko2023regularized,
doikov2024super,
han2025low,
ning2023multi,
wang2022newton,
robinson2009generalized,
mordukhovich2023globally,
wang2026damage,
robinson2009local,
shuo2026lecture,lin2023monotone,yu2026optimization,liang2025squared,jongen1987inertia,guddat1990parametric,bellon2024time,tang2022running,liu2025new,josephy1979newton,si2024riemannian,yu2026optimization1,longman2023method,yao2023relative,han2024continuous,ha1987application,yu2026pattern,kojima2009continuous,seguin2022continuation,liu1995perturbation,pang1996piecewise,liu2022iterative,eikenbroek2022improving,ralph1995directional,yu2026optimization2,scheel2000mathematical,bonnans1992developpement,bonnans1992expansion,dempe1993directional,liu2023iterative,shapiro1988sensitivity,pang1990newton,robinson1991implicit,zheng2025enhanced,hang2025smoothness,hang2024role,scholtes2012introduction,cui2022nonconvex}. 

\bibliography{reference1} 

@article{falk1995bilevel,
  title={On bilevel programming, Part I: general nonlinear cases},
  author={Falk, James E and Liu, Jiming},
  journal={Mathematical Programming},
  volume={70},
  number={1},
  pages={47--72},
  year={1995},
  publisher={Springer}
}

@article{aghasi2025fully,
  title={Fully zeroth-order bilevel programming via Gaussian smoothing},
  author={Aghasi, Alireza and Ghadimi, Saeed},
  journal={Journal of Optimization Theory and Applications},
  volume={205},
  number={2},
  pages={31},
  year={2025},
  publisher={Springer}
}

@article{hong2023two,
  title={A two-timescale stochastic algorithm framework for bilevel optimization: Complexity analysis and application to actor-critic},
  author={Hong, Mingyi and Wai, Hoi-To and Wang, Zhaoran and Yang, Zhuoran},
  journal={SIAM Journal on Optimization},
  volume={33},
  number={1},
  pages={147--180},
  year={2023},
  publisher={SIAM}
}

@article{kovccvara1994optimization,
  title={On optimization of systems governed by implicit complementarity problems},
  author={Ko{\v{c}}cvara, Michal and Outrata, Jan V},
  journal={Numerical Functional Analysis and Optimization},
  volume={15},
  number={7-8},
  pages={869--887},
  year={1994},
  publisher={Taylor \& Francis}
}

@article{chaudet2020shape,
  title={Shape derivatives for the penalty formulation of elastic contact problems with Tresca friction},
  author={Chaudet-Dumas, Bastien and Deteix, Jean},
  journal={SIAM Journal on Control and Optimization},
  volume={58},
  number={6},
  pages={3237--3261},
  year={2020},
  publisher={SIAM}
}

@article{outrata1995numerical,
  title={A numerical approach to optimization problems with variational inequality constraints},
  author={Outrata, Ji{\v{r}}{\'\i} and Zowe, Jochem},
  journal={Mathematical Programming},
  volume={68},
  number={1},
  pages={105--130},
  year={1995},
  publisher={Springer}
}

@article{cui2023complexity,
  title={Complexity guarantees for an implicit smoothing-enabled method for stochastic MPECs},
  author={Cui, Shisheng and Shanbhag, Uday V and Yousefian, Farzad},
  journal={Mathematical Programming},
  volume={198},
  number={2},
  pages={1153--1225},
  year={2023},
  publisher={Springer}
}

@article{christof2020nonsmooth,
  title={A nonsmooth trust-region method for locally Lipschitz functions with application to optimization problems constrained by variational inequalities},
  author={Christof, Constantin and De los Reyes, Juan Carlos and Meyer, Christian},
  journal={SIAM Journal on Optimization},
  volume={30},
  number={3},
  pages={2163--2196},
  year={2020},
  publisher={SIAM}
}

@article{rawat2026augmented,
  title={Augmented Lagrangian Neural Network for Solving Mathematical Programs with Equilibrium Constraints},
  author={Rawat, Anjali and Singh, Vinay},
  journal={Journal of Optimization Theory and Applications},
  volume={209},
  number={1},
  pages={15},
  year={2026},
  publisher={Springer}
}

@article{robinson1980strongly,
  title={Strongly regular generalized equations},
  author={Robinson, Stephen M},
  journal={Mathematics of Operations Research},
  volume={5},
  number={1},
  pages={43--62},
  year={1980},
  publisher={INFORMS}
}

@article{shin2023near,
  title={Near-optimal distributed linear-quadratic regulator for networked systems},
  author={Shin, Sungho and Lin, Yiheng and Qu, Guannan and Wierman, Adam and Anitescu, Mihai},
  journal={SIAM Journal on Control and Optimization},
  volume={61},
  number={3},
  pages={1113--1135},
  year={2023},
  publisher={SIAM}
}

@article{bolte2024differentiating,
  title={Differentiating nonsmooth solutions to parametric monotone inclusion problems},
  author={Bolte, J{\'e}r{\^o}me and Pauwels, Edouard and Silveti-Falls, Antonio},
  journal={SIAM Journal on Optimization},
  volume={34},
  number={1},
  pages={71--97},
  year={2024},
  publisher={SIAM}
}

@article{nghia2025geometric,
  title={Geometric characterizations of Lipschitz stability for convex optimization problems},
  author={Nghia, Tran TA},
  journal={SIAM Journal on Optimization},
  volume={35},
  number={2},
  pages={927--958},
  year={2025},
  publisher={SIAM}
}

@article{chen2025aubin,
  title={Aubin property and strong regularity are equivalent for nonlinear second-order cone programming},
  author={Chen, Liang and Chen, Ruoning and Sun, Defeng and Zhu, Junyuan},
  journal={SIAM Journal on Optimization},
  volume={35},
  number={2},
  pages={712--738},
  year={2025},
  publisher={SIAM}
}

@incollection{kojima1980strongly,
  title={Strongly stable stationary solutions in nonlinear programs},
  author={Kojima, Masakazu},
  booktitle={Analysis and computation of fixed points},
  pages={93--138},
  year={1980},
  publisher={Elsevier}
}

@article{chen2026characterizations,
  title={Characterizations of the Aubin property of the solution mapping for nonlinear semidefinite programming: L. Chen et al.},
  author={Chen, Liang and Chen, Ruoning and Sun, Defeng and Zhang, Liping},
  journal={Mathematical Programming},
  volume={215},
  number={1},
  pages={637--668},
  year={2026},
  publisher={Springer}
}

@article{cui2026lipschitz,
  title={Lipschitz Stability of Least-Squares Problems Regularized by Functions with C 2-Cone Reducible Conjugates},
  author={Cui, Ying and Hoheisel, Tim and Nghia, Tran TA and Sun, Defeng},
  journal={Mathematics of Operations Research},
  year={2026},
  publisher={INFORMS}
}

@article{shin2022exponential,
  title={Exponential decay of sensitivity in graph-structured nonlinear programs},
  author={Shin, Sungho and Anitescu, Mihai and Zavala, Victor M},
  journal={SIAM Journal on Optimization},
  volume={32},
  number={2},
  pages={1156--1183},
  year={2022},
  publisher={SIAM}
}

@misc{bank1982d,
  title={D. Klatte, B. Kummer, and K. Tammer: Non-Linear Parametric Optimization},
  author={Bank, B and Guddat, J},
  journal={Berlin},
  year={1982}
}

@article{bonnans1994local,
  title={Local analysis of Newton-type methods for variational inequalities and nonlinear programming},
  author={Bonnans, J Fr{\'e}d{\'e}ric},
  journal={Applied Mathematics and Optimization},
  volume={29},
  number={2},
  pages={161--186},
  year={1994},
  publisher={Springer}
}

@article{chen2025two,
  title={Two Typical Implementable Semismooth* Newton Methods for Generalized Equations are G-Semismooth Newton Methods},
  author={Chen, Liang and Sun, Defeng and Zhang, Wangyongquan},
  journal={Mathematics of Operations Research},
  year={2025},
  publisher={INFORMS}
}

@article{mohammadi2022variational,
  title={Variational analysis of composite models with applications to continuous optimization},
  author={Mohammadi, Ashkan and Mordukhovich, Boris S and Sarabi, M Ebrahim},
  journal={Mathematics of Operations Research},
  volume={47},
  number={1},
  pages={397--426},
  year={2022},
  publisher={INFORMS}
}

@article{dussault2026polyhedral,
  title={Polyhedral Newton-min algorithms for complementarity problems: J.-P. Dussault et al.},
  author={Dussault, Jean-Pierre and Frappier, Mathieu and Gilbert, Jean Charles},
  journal={Mathematical Programming},
  volume={215},
  number={1},
  pages={269--324},
  year={2026},
  publisher={Springer}
}

@article{qi2000constant,
  title={On the constant positive linear dependence condition and its application to SQP methods},
  author={Qi, Liqun and Wei, Zengxin},
  journal={SIAM Journal on Optimization},
  volume={10},
  number={4},
  pages={963--981},
  year={2000},
  publisher={SIAM}
}

@article{facchinei1998accurate,
  title={On the accurate identification of active constraints},
  author={Facchinei, Francisco and Fischer, Andreas and Kanzow, Christian},
  journal={SIAM Journal on Optimization},
  volume={9},
  number={1},
  pages={14--32},
  year={1998},
  publisher={SIAM}
}

@incollection{jittorntrum2009solution,
  title={Solution point differentiability without strict complementarity in nonlinear programming},
  author={Jittorntrum, Krisorn},
  booktitle={Sensitivity, Stability and Parametric Analysis},
  pages={127--138},
  year={2009},
  publisher={Springer}
}

@article{de2023function,
  title={A function approximation approach for parametric optimization},
  author={De Marchi, Alberto and Dreves, Axel and Gerdts, Matthias and Gottschalk, Simon and Rogovs, Sergejs},
  journal={Journal of Optimization Theory and Applications},
  volume={196},
  number={1},
  pages={56--77},
  year={2023},
  publisher={Springer}
}

@article{kyparisis1990sensitivity,
  title={Sensitivity analysis for nonlinear programs and variational inequalities with nonunique multipliers},
  author={Kyparisis, Jerzy},
  journal={Mathematics of operations research},
  volume={15},
  number={2},
  pages={286--298},
  year={1990},
  publisher={INFORMS}
}

@article{kyparisis1992parametric,
  title={Parametric variational inequalities with multivalued solution sets},
  author={Kyparisis, Jerzy},
  journal={Mathematics of Operations Research},
  volume={17},
  number={2},
  pages={341--364},
  year={1992},
  publisher={INFORMS}
}

@article{liu1995sensitivity,
  title={Sensitivity analysis in nonlinear programs and variational inequalities via continuous selections},
  author={Liu, Jiming},
  journal={SIAM Journal on Control and Optimization},
  volume={33},
  number={4},
  pages={1040--1060},
  year={1995},
  publisher={SIAM}
}

@article{pang11995stability,
  title={Stability of Parametric Nonsmooth Equations},
  author={Pang$^1$, Jong-Shi},
  journal={Recent Advances in Nonsmooth Optimization},
  pages={261},
  year={1995},
  publisher={World Scientific}
}

@article{qiu1992sensitivity,
  title={Sensitivity analysis for variational inequalities},
  author={Qiu, Yuping and Magnanti, Thomas L},
  journal={Mathematics of Operations Research},
  volume={17},
  number={1},
  pages={61--76},
  year={1992},
  publisher={INFORMS}
}

@article{aussel2024variational,
  title={Variational and quasi-variational inequalities under local reproducibility: solution concept and applications},
  author={Aussel, Didier and Chaipunya, Parin},
  journal={Journal of Optimization Theory and Applications},
  volume={203},
  number={2},
  pages={1531--1563},
  year={2024},
  publisher={Springer}
}

@article{reinoza1985strong,
  title={The strong positivity conditions},
  author={Reinoza, Alfonso},
  journal={Mathematics of operations research},
  volume={10},
  number={1},
  pages={54--62},
  year={1985},
  publisher={INFORMS}
}

@book{facchinei2003finite,
  title={Finite-dimensional variational inequalities and complementarity problems},
  author={Facchinei, Francisco and Pang, Jong-Shi},
  year={2003},
  publisher={Springer}
}

@article{harker1990finite,
  title={Finite-dimensional variational inequality and nonlinear complementarity problems: a survey of theory, algorithms and applications},
  author={Harker, Patrick T and Pang, Jong-Shi},
  journal={Mathematical programming},
  volume={48},
  number={1},
  pages={161--220},
  year={1990},
  publisher={Springer}
}

@misc{kleinmichel1972av,
  title={AV Fiacco/GP McCormick, Nonlinear Programming: Sequential Unconstrained Minimization Techniques. XIV+ 210 S. m. Fig. New York/London/Sydney/Toronto 1969. John Wiley and Sons, Inc. Preis. geb.{\pounds} 5, 20},
  author={Kleinmichel, H},
  year={1972},
  publisher={Wiley Online Library}
}

@book{ortega2000iterative,
  title={Iterative solution of nonlinear equations in several variables},
  author={Ortega, James M and Rheinboldt, Werner C},
  year={2000},
  publisher={SIAM}
}

@book{bai2021matrix,
  title={Matrix analysis and computations},
  author={Bai, Zhong-Zhi and Pan, Jian-Yu},
  year={2021},
  publisher={SIAM}
}

@article{gander2026landmarks,
  title={Landmarks in the history of iterative methods},
  author={Gander, Martin J and Henry, Philippe and Wanner, Gerhard},
  journal={SIAM Review},
  volume={68},
  number={1},
  pages={3--90},
  year={2026},
  publisher={SIAM}
}

@article{mishchenko2023regularized,
  title={Regularized Newton method with global convergence},
  author={Mishchenko, Konstantin},
  journal={SIAM Journal on Optimization},
  volume={33},
  number={3},
  pages={1440--1462},
  year={2023},
  publisher={SIAM}
}

@article{doikov2024super,
  title={Super-universal regularized Newton method},
  author={Doikov, Nikita and Mishchenko, Konstantin and Nesterov, Yurii},
  journal={SIAM Journal on Optimization},
  volume={34},
  number={1},
  pages={27--56},
  year={2024},
  publisher={SIAM}
}

@article{han2025low,
  title={A low-rank admm splitting approach for semidefinite programming},
  author={Han, Qiushi and Li, Chenxi and Lin, Zhenwei and Chen, Caihua and Deng, Qi and Ge, Dongdong and Liu, Huikang and Ye, Yinyu},
  journal={INFORMS Journal on Computing},
  year={2025},
  publisher={INFORMS}
}

@incollection{robinson2009generalized,
  title={Generalized equations and their solutions, Part I: Basic theory},
  author={Robinson, Stephen M},
  booktitle={Point-to-Set Maps and Mathematical Programming},
  pages={128--141},
  year={2009},
  publisher={Springer}
}

@article{ning2023multi,
  title={Multi-dimensional path-dependent forward-backward stochastic variational inequalities},
  author={Ning, Ning and Wu, Jing},
  journal={Set-Valued and Variational Analysis},
  volume={31},
  number={1},
  pages={2},
  year={2023},
  publisher={Springer}
}

@article{mordukhovich2023globally,
  title={A globally convergent proximal Newton-type method in nonsmooth convex optimization},
  author={Mordukhovich, Boris S and Yuan, Xiaoming and Zeng, Shangzhi and Zhang, Jin},
  journal={Mathematical Programming},
  volume={198},
  number={1},
  pages={899--936},
  year={2023},
  publisher={Springer}
}

@article{wang2022newton,
  title={Newton’s method for solving generalized equations without Lipschitz condition},
  author={Wang, Jiaxi and Ouyang, Wei},
  journal={Journal of Optimization Theory and Applications},
  volume={192},
  number={2},
  pages={510--532},
  year={2022},
  publisher={Springer}
}

@incollection{robinson2009local,
  title={Local structure of feasible sets in nonlinear programming, Part III: Stability and sensitivity},
  author={Robinson, Stephen M},
  booktitle={Nonlinear Analysis and Optimization},
  pages={45--66},
  year={2009},
  publisher={Springer}
}

@article{lin2023monotone,
  title={Monotone inclusions, acceleration, and closed-loop control},
  author={Lin, Tianyi and Jordan, Michael I},
  journal={Mathematics of Operations Research},
  volume={48},
  number={4},
  pages={2353--2382},
  year={2023},
  publisher={INFORMS}
}

@article{liang2025squared,
  title={A squared smoothing Newton method for semidefinite programming},
  author={Liang, Ling and Sun, Defeng and Toh, Kim-Chuan},
  journal={Mathematics of Operations Research},
  volume={50},
  number={4},
  pages={2873--2908},
  year={2025},
  publisher={INFORMS}
}

@article{jongen1987inertia,
  title={On inertia and Schur complement in optimization},
  author={Jongen, H Th and M{\"o}bert, T and R{\"u}ckmann, J and Tammer, K},
  journal={Linear Algebra and its Applications},
  volume={95},
  pages={97--109},
  year={1987},
  publisher={Elsevier}
}

@book{guddat1990parametric,
  title={Parametric optimization: singularities, pathfollowing and jumps},
  author={Guddat, J{\"u}rgen and Vazquez, F Guerra and Jongen, Hubertus Th},
  year={1990},
  publisher={Springer}
}

@article{bellon2024time,
  title={Time-varying semidefinite programming: Path following a Burer--Monteiro factorization},
  author={Bellon, Antonio and Dressler, Mareike and Kungurtsev, Vyacheslav and Mare{\v{c}}ek, Jakub and Uschmajew, Andr{\'e}},
  journal={SIAM Journal on Optimization},
  volume={34},
  number={1},
  pages={1--26},
  year={2024},
  publisher={SIAM}
}

@article{tang2022running,
  title={Running primal-dual gradient method for time-varying nonconvex problems},
  author={Tang, Yujie and Dall'Anese, Emiliano and Bernstein, Andrey and Low, Steven},
  journal={SIAM Journal on Control and Optimization},
  volume={60},
  number={4},
  pages={1970--1990},
  year={2022},
  publisher={SIAM}
}

@article{liu2025new,
  title={New methods for parametric optimization via differential equations},
  author={Liu, Heyuan and Grigas, Paul},
  journal={SIAM Journal on Optimization},
  volume={35},
  number={3},
  pages={1524--1550},
  year={2025},
  publisher={SIAM}
}

@techreport{josephy1979newton,
  title={Newton's Method for Generalized Equations.},
  author={Josephy, Norman H},
  year={1979}
}

@article{si2024riemannian,
  title={A Riemannian proximal Newton method},
  author={Si, Wutao and Absil, P-A and Huang, Wen and Jiang, Rujun and Vary, Simon},
  journal={SIAM Journal on Optimization},
  volume={34},
  number={1},
  pages={654--681},
  year={2024},
  publisher={SIAM}
}

@article{gowda1994stability,
  title={Stability analysis of variational inequalities and nonlinear complementarity problems, via the mixed linear complementarity problem and degree theory},
  author={Gowda, M Seetharama and Pang, Jong-Shi},
  journal={Mathematics of Operations Research},
  volume={19},
  number={4},
  pages={831--879},
  year={1994},
  publisher={INFORMS}
}

@article{yao2023relative,
  title={Relative Lipschitz-like property of parametric systems via projectional coderivatives},
  author={Yao, Wenfang and Yang, Xiaoqi},
  journal={SIAM Journal on Optimization},
  volume={33},
  number={3},
  pages={2021--2040},
  year={2023},
  publisher={SIAM}
}

@article{han2024continuous,
  title={Continuous selections of solutions to parametric variational inequalities},
  author={Han, Shaoning and Pang, Jong-Shi},
  journal={SIAM Journal on Optimization},
  volume={34},
  number={1},
  pages={870--892},
  year={2024},
  publisher={SIAM}
}

@article{ha1987application,
  title={Application of degree theory in stability of the complementarity problem},
  author={Ha, Cu Duong},
  journal={Mathematics of Operations Research},
  volume={12},
  number={2},
  pages={368--376},
  year={1987},
  publisher={INFORMS}
}

@incollection{kojima2009continuous,
  title={Continuous deformation of nonlinear programs},
  author={Kojima, Masakazu and Hirabayashi, Ryuichi},
  booktitle={Sensitivity, Stability and Parametric Analysis},
  pages={150--198},
  year={2009},
  publisher={Springer}
}

@article{seguin2022continuation,
  title={Continuation methods for Riemannian optimization},
  author={S{\'e}guin, Axel and Kressner, Daniel},
  journal={SIAM Journal on Optimization},
  volume={32},
  number={2},
  pages={1069--1093},
  year={2022},
  publisher={SIAM}
}

@phdthesis{liu1995perturbation,
  title={Perturbation analysis in nonlinear programs and variational inequalities},
  author={Liu, Jiming},
  year={1995},
  school={George Washington University}
}

@article{pang1996piecewise,
  title={Piecewise smoothness, local invertibility, and parametric analysis of normal maps},
  author={Pang, Jong-Shi and Ralph, Daniel},
  journal={Mathematics of operations research},
  volume={21},
  number={2},
  pages={401--426},
  year={1996},
  publisher={INFORMS}
}

@article{eikenbroek2022improving,
  title={Improving the performance of a traffic system by fair rerouting of travelers},
  author={Eikenbroek, Oskar AL and Still, Georg J and Van Berkum, Eric C},
  journal={European Journal of Operational Research},
  volume={299},
  number={1},
  pages={195--207},
  year={2022},
  publisher={Elsevier}
}

@article{ralph1995directional,
  title={Directional derivatives of the solution of a parametric nonlinear program},
  author={Ralph, Daniel and Dempe, Stephan},
  journal={Mathematical programming},
  volume={70},
  number={1},
  pages={159--172},
  year={1995},
  publisher={Springer}
}

@article{scheel2000mathematical,
  title={Mathematical programs with complementarity constraints: Stationarity, optimality, and sensitivity},
  author={Scheel, Holger and Scholtes, Stefan},
  journal={Mathematics of Operations Research},
  volume={25},
  number={1},
  pages={1--22},
  year={2000},
  publisher={INFORMS}
}

@article{bonnans1992developpement,
  title={D{\'e}veloppement de solutions exactes et approch{\'e}es en programmation nonlin{\'e}aire},
  author={Bonnans, JF and Ioffe, AD and Shapiro, Alexander},
  journal={Comptes rendus de l'Acad{\'e}mie des sciences. S{\'e}rie 1, Math{\'e}matique},
  volume={315},
  number={2},
  pages={119--123},
  year={1992}
}

@inproceedings{bonnans1992expansion,
  title={Expansion of exact and approximate solutions in nonlinear programming},
  author={Bonnans, J Fr{\'e}d{\'e}ric and Ioffe, Alexander D and Shapiro, Alexander},
  booktitle={Advances in Optimization: Proceedings of the 6th French-German Colloquium on Optimization Held at Lambrecht, FRG, June 2--8, 1991},
  pages={103--117},
  year={1992},
  organization={Springer}
}

@article{dempe1993directional,
  title={Directional differentiability of optimal solutions under Slater's condition},
  author={Dempe, Stephan},
  journal={Mathematical programming},
  volume={59},
  number={1},
  pages={49--69},
  year={1993},
  publisher={Springer}
}

@article{shapiro1988sensitivity,
  title={Sensitivity analysis of nonlinear programs and differentiability properties of metric projections},
  author={Shapiro, Alexander},
  journal={SIAM Journal on Control and Optimization},
  volume={26},
  number={3},
  pages={628--645},
  year={1988},
  publisher={SIAM}
}

@article{pang1990newton,
  title={Newton's method for B-differentiable equations},
  author={Pang, Jong-Shi},
  journal={Mathematics of operations research},
  volume={15},
  number={2},
  pages={311--341},
  year={1990},
  publisher={INFORMS}
}

@article{khanh2024globally,
  title={Globally convergent coderivative-based generalized Newton methods in nonsmooth optimization},
  author={Khanh, Pham Duy and Mordukhovich, Boris S and Phat, Vo Thanh and Tran, Dat Ba},
  journal={Mathematical Programming},
  volume={205},
  number={1},
  pages={373--429},
  year={2024},
  publisher={Springer}
}

@article{robinson1991implicit,
  title={An implicit-function theorem for a class of nonsmooth functions},
  author={Robinson, Stephen M},
  journal={Mathematics of operations research},
  volume={16},
  number={2},
  pages={292--309},
  year={1991},
  publisher={INFORMS}
}

@article{hang2025smoothness,
  title={Smoothness of subgradient mappings and its applications in parametric optimization},
  author={Hang, Nguyen TV and Sarabi, Ebrahim},
  journal={Set-Valued and Variational Analysis},
  volume={33},
  number={4},
  pages={41},
  year={2025},
  publisher={Springer}
}

@article{hang2024role,
  title={Role of subgradients in variational analysis of polyhedral functions},
  author={Hang, Nguyen TV and Jung, Woosuk and Sarabi, Ebrahim},
  journal={Journal of Optimization Theory and Applications},
  volume={200},
  number={3},
  pages={1160--1192},
  year={2024},
  publisher={Springer}
}

@book{scholtes2012introduction,
  title={Introduction to piecewise differentiable equations},
  author={Scholtes, Stefan},
  year={2012},
  publisher={Springer}
}

@article{cui2022nonconvex,
  title={Nonconvex and nonsmooth approaches for affine chance-constrained stochastic programs},
  author={Cui, Ying and Liu, Junyi and Pang, Jong-Shi},
  journal={Set-Valued and Variational Analysis},
  volume={30},
  number={3},
  pages={1149--1211},
  year={2022},
  publisher={Springer}
}

@article{liu2025bidirectional,
  title={Bidirectional endothelial feedback drives turing-vascular patterning and drug-resistance niches: a hybrid PDE-agent-based study},
  author={Liu, Zonghao and Wang, Louis Shuo and Yu, Jiguang and Zhang, Jilin and Martel, Erica and Li, Shijia},
  journal={Bioengineering},
  volume={12},
  number={10},
  pages={1097},
  year={2025},
  publisher={MDPI}
}

@article{wang2025analysis,
  title={Analysis framework for stochastic predator--prey model with demographic noise},
  author={Wang, Louis Shuo and Yu, Jiguang},
  journal={Results in Applied Mathematics},
  volume={27},
  pages={100621},
  year={2025},
  publisher={Elsevier}
}

@article{wang2026algebraic,
  title={Algebraic--spectral thresholds and discrete--continuous stability transfer in Leslie--Gower systems},
  author={Wang, Louis Shuo and Yu, Jiguang},
  journal={Electronic Research Archive},
  volume={34},
  number={1},
  pages={251--290},
  year={2026}
}

@article{wang2025analysis1,
  title={Analysis and Mean-Field Limit of a Hybrid PDE-ABM Modeling Angiogenesis-Regulated Resistance Evolution},
  author={Wang, Louis Shuo and Yu, Jiguang and Li, Shijia and Liu, Zonghao},
  journal={Mathematics},
  volume={13},
  number={17},
  pages={2898},
  year={2025},
  publisher={MDPI}
}

@article{wang2026damage,
  title={A damage-structured PDE model of stem cell hierarchies: The dual role of dedifferentiation in tissue homeostasis and aging},
  author={Wang, Louis Shuo and Yu, Jiguang and Liu, Zonghao},
  journal={Plos one},
  volume={21},
  number={2},
  pages={e0335163},
  year={2026},
  publisher={Public Library of Science San Francisco, CA USA}
}

@article{yu2026pattern,
  title={Pattern suppression and recovery under one-way versus two-way chemotactic coupling in hybrid partial differential equation--ordinary differential equation models},
  author={Yu, Jiguang and Wang, Louis Shuo and Liu, Zonghao and Liu, Jingfeng},
  journal={Transport Phenomena},
  number={0},
  year={2026},
  publisher={De Gruyter}
}

@article{shuo2026lecture,
  title={Lecture Note for Bounded Controls in Continuous-Time and Control of Several Variables},
  author={Shuo Wang, Louis},
  journal={arXiv e-prints},
  pages={arXiv--2604},
  year={2026}
}

@inproceedings{liu2023auxiliary,
  title={Auxiliary-variable adaptive control barrier functions for safety critical systems},
  author={Liu, Shuo and Xiao, Wei and Belta, Calin A},
  booktitle={2023 62nd IEEE Conference on Decision and Control (CDC)},
  pages={8602--8607},
  year={2023},
  organization={IEEE}
}

@inproceedings{chen2025control,
  title={Control barrier functions via minkowski operations for safe navigation among polytopic sets},
  author={Chen, Yi-Hsuan and Liu, Shuo and Xiao, Wei and Belta, Calin and Otte, Michael},
  booktitle={2025 IEEE 64th Conference on Decision and Control (CDC)},
  pages={4481--4488},
  year={2025},
  organization={IEEE}
}

@inproceedings{liu2024auxiliary,
  title={Auxiliary-variable adaptive control lyapunov barrier functions for spatio-temporally constrained safety-critical applications},
  author={Liu, Shuo and Xiao, Wei and Belta, Calin A},
  booktitle={2024 IEEE 63rd Conference on Decision and Control (CDC)},
  pages={8098--8104},
  year={2024},
  organization={IEEE}
}

@article{liu2025learning,
  title={Learning-Enabled Iterative Convex Optimization for Safety-Critical Model Predictive Control},
  author={Liu, Shuo and Huang, Zhe and Zeng, Jun and Sreenath, Koushil and Belta, Calin A},
  journal={IEEE Open Journal of Control Systems},
  year={2025},
  publisher={IEEE}
}

@inproceedings{liu2024feasibility,
  title={Feasibility-guaranteed safety-critical control with applications to heterogeneous platoons},
  author={Liu, Shuo and Xiao, Wei and Belta, Calin},
  booktitle={2024 IEEE 63rd Conference on Decision and Control (CDC)},
  pages={8066--8073},
  year={2024},
  organization={IEEE}
}

@article{lin2026hierarchical,
  title={Hierarchical Multi-Agent MCTS for Safety-Critical Coordination in Mixed-Autonomy Roundabouts},
  author={Lin, Zhihao and Lan, Jianglin and Liu, Shuo and Tian, Zhen and Zhao, Dezong and Wei, Chongfeng},
  journal={IEEE Transactions on Vehicular Technology},
  year={2026},
  publisher={IEEE}
}

@inproceedings{liu2025risk,
  title={Risk-aware adaptive control barrier functions for safe control of nonlinear systems under stochastic uncertainty},
  author={Liu, Shuo and Belta, Calin A},
  booktitle={2025 IEEE 64th Conference on Decision and Control (CDC)},
  pages={4875--4881},
  year={2025},
  organization={IEEE}
}

@article{longman2023method,
  title={A method to speed up convergence of iterative learning control for high precision repetitive motions},
  author={Longman, Richard W and Liu, Shuo and Elsharhawy, Tarek A},
  journal={arXiv preprint arXiv:2307.15912},
  year={2023}
}

@article{liu2022iterative,
  title={Iterative convex optimization for model predictive control with discrete-time high-order control barrier functions},
  author={Liu, Shuo and Zeng, Jun and Sreenath, Koushil and Belta, Calin A},
  journal={arXiv preprint arXiv:2210.04361},
  year={2022}
}

@article{zheng2025enhanced,
  title={Enhanced mean field game for interactive decision-making with varied stylish multi-vehicles},
  author={Zheng, Liancheng and Tian, Zhen and He, Yangfan and Liu, Shuo and Chen, Huilin and Yuan, Fujiang and Peng, Yanhong},
  journal={arXiv preprint arXiv:2509.00981},
  year={2025}
}

@article{wang2025multi,
  title={Multi-strategy Hybrid Improved Intelligent Algorithm for Solving UAV-MTSP},
  author={Wang, Zixin and Wang, Danqing and Yu, Jiguang},
  journal={Information Technology and Control},
  volume={54},
  number={2},
  pages={413--438},
  year={2025}
}

@inproceedings{gao2022rolling,
  title={Rolling prediction model of closing price based on EEMD data noise reduction and HGS-DELM},
  author={Gao, Yuansheng and Li, Lei and Yu, Jiguang},
  booktitle={2022 International Conference on Data Analytics, Computing and Artificial Intelligence (ICDACAI)},
  pages={255--260},
  year={2022},
  organization={IEEE}
}

@article{yu2026optimization,
  title={Optimization Workshop Notes for Mathematical Programming with Equilibrium Constraints (MPECs): Second-Order Optimality Conditions},
  author={Yu, Jiguang},
  journal={arXiv preprint arXiv:2604.20992},
  year={2026}
}

@article{yu2026optimization1,
  title={Optimization Workshop Notes for Mathematical Programming with Equilibrium Constraints (MPECs): Verification of MPEC Hypotheses},
  author={Yu, Jiguang},
  journal={arXiv preprint arXiv:2604.20988},
  year={2026}
}

@article{yu2026optimization2,
  title={Optimization Workshop Notes for Mathematical Programming with Equilibrium Constraints Algorithms: Penalty Interior-Point, Implicit-Programming, and Piecewise SQP},
  author={Yu, Jiguang},
  journal={arXiv preprint arXiv:2604.15690},
  year={2026}
}

@inproceedings{liu2023iterative,
  title={Iterative Convex Optimization for Model Predictive Control with Discrete-Time High-Order Control Barrier Functions},
  author={Liu, Shuo and Zeng, Jun and Sreenath, Koushil and Belta, Calin A},
  booktitle={American Control Conference},
  year={2023},
  organization={IEEE}
}

\end{document}